\documentclass[12pt]{article}
\usepackage{amsfonts}
\usepackage{amsthm,amsmath,amssymb,anysize}
\newtheorem{lemma}{Lemma}[section]
\newtheorem{theorem}[lemma]{Theorem}

\newtheorem{example}[lemma]{Example}
\marginsize{3.1cm}{3.1cm}{3.1cm}{3cm}
\numberwithin{equation}{section}

\newcommand {\bcdot}   {\mathbin{\hbox{\raise.4ex\hbox{\bf.}}}} 

\def \opname#1#2%
  {\expandafter\newcommand \csname #1\endcsname {\mathop{\mathrm{#2}}\nolimits }}

\opname{Span} {{Span}}

\title{\textsf{ Cohomology of  Heisenberg Lie superalgebras}}

\author{\textsc{Wei BAI}\; {and} \textsc{Wende LIU}
\footnote{Corresponding author. Email: wendeliu@ustc.edu.cn} \setcounter{footnote}{-1}
\footnote{Supported by  NSF of  HLJP (A2015003, A201412) and NSF of China (11501151, 11471090)
}\\
  \\
\ \ \ {\small\textit{School of Mathematical Sciences},
\textit{Harbin Normal University}} \\
{\small\textit{Harbin 150025, China}}
 }

\begin{document}
\maketitle
{\small\begin{quotation}
\noindent\textbf{Abstract} Suppose the ground field to be algebraically closed and of  characteristic different from  $2$ and $3$. All Heisenberg Lie superalgebras consist of  two super versions of the
Heisenberg Lie algebras, $\frak{h}_{2m,n}$ and $\frak{ba}_n$ with $m$  a nonnegative integer and $n$ a positive integer. The space of a ``classical" Heisenberg Lie
superalgebra $\frak{h}_{2m,n}$ is the direct sum of a superspace with a non-degenerate
anti-supersymmetric even bilinear form and a one-dimensional space of
values of this form constituting the even center. The other super analog
of the Heisenberg Lie algebra, $\frak{ba}_n$, is constructed by means of a
non-degenerate anti-supersymmetric odd bilinear form with values in the one-dimensional
odd center. In this paper, we study the  cohomology  of $\frak{h}_{2m,n}$ and $\frak{ba}_n$ with coefficients in the trivial
module by using  the Hochschild-Serre spectral sequences  relative to a suitable ideal.  In  characteristic zero case, for any Heisenberg Lie superalgebra, we determine completely  the Betti numbers and associative superalgebra structure for
 their cohomology. In  characteristic $p>3$ case, we determine  the associative superalgebra structures   for
 the divided power cohomology  of  $\frak{ba}_n$ and we also make an attempt to determine  the cohomology of $\frak{h}_{2m,n}$ by computing it in a low-dimensional case.
\\


  \end{quotation}}

  \setcounter{section}{-1}

\section{Introduction}\label{nco-h000}

The theory of Lie superalgebras and Lie supergroups has many applications in various areas of physics and cohomology is an important tool in the modern
mathematics and theoretical physics (for example, see \cite{GM, Musson, WLS}).  One may find many studies on Lie superalgebra cohomology and its application in the literature (for example, see  \cite{G1,G2,GL,Leb1,LLS}).

For a given Lie superalgebra, generally speaking, it is not easy to determine the  Betti numbers,  superalgebra structures of cohomology, and so on. We should mention that  the cohomology space with coefficients in the trivial
module for a Lie superalgebra possesses a natural associative superalgebra structure arising from the cup product of the corresponding cochain space
\cite{Leites, Musson}. A standard fact is that such an associative superalgebra is graded-supercommutative  (for a definition, see Section 1) for a Lie superalgebra with nontrivial odd part and is supercommutative for a Lie algebra. The present paper aims  to determine explicitly the Betti numbers and associative superalgebra structures of the cohomology spaces with coefficients in the trivial module for a class of nilpotent Lie superalgebras---Heisenberg Lie superalgebras, over an algebraically closed field of characteristic different from $2$ and $3.$

Not as in the super case, for Heisenberg Lie algebras,  one has obtained many conclusions on cohomology
with coefficients in the trivial module.
Santharoubane \cite{LJS} studied  the cohomology  of Heisenberg Lie
algebras   and gave a description of cocycles, coboundaries and
cohomology spaces  over a field of characteristic zero.
Sk\"oldberg \cite{Emil} calculated   Betti numbers, which are
dimensions of $n$-cohomology spaces,  for Heisenberg Lie algebras
over  a field  of characteristic two. Cairns and Jambor \cite{CJ}
extended Sk\"oldberg's work \cite{Emil}   to the case of
arbitrary prime characteristic.

\section{Heisenberg Lie Superalgebras}\label{nco-h001}

Throughout this paper, the ground  field $\mathbb{F}$ is algebraically closed and of  characteristic $p\not=2,3$.
A Heisenberg Lie superalgebra is by definition a two-step nilpotent Lie superalgebra
with 1-dimensional center  \cite{RSS}.
Note that, over a non-algebraically closed field, there are several equivalence
classes of non-degenerate symmetric bilinear forms on a given space.
 To simplify our life, we assume in this paper the ground field $\mathbb{F}$ is
algebraically closed. In this case, all finite dimensional
Heisenberg Lie superalgebras split precisely into two types
$\frak{h}_{2m,n}$ and $\frak{ba}_{n}$, where  $0\leq m\in\mathbb{Z}$ and $1\leq n\in\mathbb{Z}$, described as
follows (we write the even basis elements  first separated from the
odd ones by a vertical bar):
\begin{itemize}

\item $\frak{h}_{2m,n}$ having a basis $
\{z; x_1,\ldots, x_m, x_{m+1},\ldots,x_{2m}\mid y_1,\ldots, y_n\}$
with non-zero brackets of basis elements (in particular, $z$
spans the even center):
\[[x_i, x_{m+i}]=z,\, [y_j, y_j]=z \text{~~for any~$1\leq i\leq m$
and $1\leq j\leq n$}.\]
\item $\frak{ba}_{n}$ having a basis $
\{x_1,\ldots, x_n\mid y_1,\ldots, y_n; z\}$ with non-zero
brackets of basis elements  (in particular, $z$ spans the odd
center): \[[x_i, y_{i}]=z\text{~~for any~~} 1\leq i\leq n.\]
\end{itemize}

A Heisenberg Lie superalgebra can
be realized to a  subalgebra of   Lie superalgebras of (odd)  contact vector fields, which is precisely the negative part with respect to a certain natural
grading \cite{Leb2}.   We should also note that Heisenberg Lie superalgebras  correspond to Heisenberg supergroups, on which one find certain recent work (for example, see \cite{BGT, CKN, T}).

The cochain space forms a graded-supercommutative superalgebra---an analog of
the polynomial algebra. We will first describe  the natural superalgebra structures on the cohomology spaces  $H^{\bullet}(\frak{h}_{2m,n})$ and $H^{\bullet}(\frak{ba}_{n})$  in case of characteristic zero.  We know that the polynomial algebra has a divided power analog in the modular case, so does the cochain space. For Lie algebras, this phenomenon  does not
 exist, since the graded-supercommutative superalgebra structure of cochain is generated by purely odd elements. Then,  we will determine  the   superalgebra structure  on the divided cohomology space $DPH^{\bullet}(\frak{ba}_{n})$ of
$\frak{ba}_{n}$
when $p>3$. For the classical Heisenberg Lie superalgebra $\frak{h}_{2m,n}$, we  make an attempt by determining  the   superalgebra
structures  on the divided cohomology space  $DPH^{\bullet}(\frak{h}_{2m,n})$  for $m=0,1$ and $ n=1$.  A complete description $DPH^{\bullet}(\frak{h}_{2m,n})$ is beyond the scope of the present paper.

 Our study is facilitated by the Hochschild-Serre spectral sequence,  of which  limit term is $E_3$ for Heisenberg Lie superalgebras (see Theorems \ref{Heisenbergth1}, \ref{Heisenbergth2} and \ref{Heisenbergth3}, and their proofs).

\section{Cohomology and Divided Power Cohomology}\label{nco-h002}
For a detailed description of superalgebras, the reader is referred
to \cite{F,Musson}.
Let $T(V)$ be the  tensor algebra of a superspace $V=V_{\bar{0}}\bigoplus V_{\bar{1}}$. Then $T(V)$
is an associative superalgebra  with a  $\mathbb{Z}_2$-grading
  induced by the $\mathbb{Z}_2$-grading of $V$. Note that $T(V)$ has also a $\mathbb{Z}$-grading structure given
by setting $\|v\|=1$ for any $v\in V$. \textit{Hereafter  $\|x\|$  denotes the $\mathbb{Z}$-degree  of a $\mathbb{Z}$-homogeneous element $x$ in a
$\mathbb{Z}$-graded  superspace.}
Let $J(V)$ denote the two-sided ideal of $T(V)$ generated by all elements of  form
$$u\otimes v+(-1)^{|u||v|}v\otimes u,\ \ u, v\in V.$$
\textit{Hereafter   $|x|$  denotes the parity of a $\mathbb{Z}_2$-homogeneous element $x$ in  a superspace.}
Note that the ideal $J(V)$ is both $\mathbb{Z}_2$-graded and $\mathbb{Z}$-graded. By definition, the super-exterior algebra is the quotient algebra of $T(V)$ by $J(V)$, denoted by $${\bigwedge} V=T(V)/J(V).$$
Then ${\bigwedge} V$ is graded-supercommutative in the following sense:
A $\mathbb{Z}$-graded associative superalgebra $A$ is said to be \textit{graded-supercommutative,} if for any $ a,b\in A$,
\begin{align*}
ab=(-1)^{\|a\|\|b\|}(-1)^{|a||b|}ba.
\end{align*}
Note that the dimension of  $k$-homogeneous subspace of the super-exterior algebra with $r$ even
 generators and $s$ odd generators,
  denoted by $\frak{d}^k(r,s)$, is
\begin{align*}
\sum_{i=0}^k\left(\begin{array}{c}
                                        r \\
                                        k-i
                                      \end{array}
 \right)\left(\begin{array}{c}
                                        s+i-1 \\
                                        i
                                      \end{array}
 \right).
\end{align*}

Over a field of characteristic zero, cohomology of Lie superalgebras is
defined by applying the Sign Rule to the corresponding definition  for Lie
algebras (for example, see \cite{F,Leites,Musson}). Let $\frak{g}$ be a Lie superalgebra and denote by $\frak{g}^*$ the dual superspace of $\frak{g}$. Then
$\bigwedge^{\bullet} \frak{g}^* =\bigoplus_{i\in\mathbb{Z}}\bigwedge^{i}\frak{g}^*$ with the differential
$\mathrm{d}$ induced by the dual of the Lie superalgebra bracket map,
$\bigwedge^{1} \frak{g}^*\longrightarrow \bigwedge^{2} \frak{g}^*$, becomes a cochain
complex, which is isomorphic to the ordinary cochain complex. In Section 2, we  determine the cohomology of complex $(\bigwedge^{\bullet} \frak{g}^*,
\mathrm{d})$ for Heisenberg Lie superalgebras in characteristic 0 case.

Over a field of prime characteristic  $p>3$, we introduce the definition of divided power cohomology of Lie superalgebra $\frak{g}$. The reader is refereed to \cite{BGLL} for further detailed information.
Fix  two $n$-tuples of positive
integers $ \underline{t}=\left(t_{1},t_{2},\ldots,t_n\right) $ and
$\pi=\left(\pi _1,\pi _2,\ldots,\pi _n\right), $ where $\pi
_i=p^{t_i}.$ Let
$$\frak{A}(m,n;\underline{t})=\mathrm{span}_{\mathbb{F}}\{u^{(\alpha)}\mid \alpha\in \mathbb{A}(m,n;\underline{t})\}$$
where
\begin{align*}\mathbb{A}(m,n;\underline{t})=\left\{\alpha=(\alpha_1,\ldots,\alpha_{m+n})\in\mathbb{N}^{m+n}\left|\begin{array}{l}
                                                                                                                 \alpha_i=0\mbox{ or } 1, 1\leq i\leq m;\\
                                                                                                                  \alpha_{m+i}<\pi_i, 1\leq i\leq n
                                                                                                               \end{array}\right.\right\}.
\end{align*}
The space $\frak{A}(m,n;\underline{t})$ is  $\mathbb{Z}_2\times\mathbb{Z}$-graded by letting $|u^{(\alpha)}|=\alpha_{m+1}+\cdots+\alpha_{m+n}$ and
$\|u^{(\alpha)}\|=\alpha_{1}+\cdots+\alpha_{m+n}$. Then
\begin{align*}&\frak{A}(m,n;\underline{t})=\bigoplus_{i=0}^{m-n+\sum_{k=1}^np^k}\frak{A}_i(m,n;\underline{t}).
\end{align*}
Moreover, $\frak{A}(m,n;\underline{t})$ is an associative  superalgebra with multiplication
\begin{align}\label{co-1018}
u^{(\alpha)}u^{(\beta)}=\lambda(\alpha, \beta)\left(\begin{array}{c}
                           \alpha+\beta\\
                           \alpha
                         \end{array}\right)u^{(\alpha+\beta)},
\end{align}
where
\begin{align*}
&\lambda(\alpha, \beta)=\left(\prod_{i=1}^m \min(1, 2-\alpha_i-\beta_i)\right)(-1)^{\sum\limits_{j=1}^m\sum\limits_{k=m+1}^{m+n}\beta_{j}\alpha_{k}+\sum\limits_{1\leq
j<k\leq m}\beta_{j}\alpha_{k}},\\
&\left(\begin{array}{c}
                           \alpha+\beta\\
                           \alpha
                         \end{array}\right)=\prod\limits_{i=1}^{m+n}\left(\begin{array}{c}
                           \alpha_{i}+\beta_{i}\\
                           \alpha_{i}
                         \end{array}\right).\end{align*}
For $\varepsilon_i=(\delta_{i,1},\ldots,\delta_{i,m+n})$, we know that $u^{(\alpha)}=u^{(\alpha_1\varepsilon_1)}\cdots u^{(\alpha_{m+n}\varepsilon_{m+n})}$,
and write ${\xi_i}$ for  $u^{(\varepsilon_i)}$, $1\leq i\leq m+n$. We call $\frak{A}(m,n;\underline{t})$ \textit{the super-divided power algebra} with
variables $\{\xi_1,\ldots,\xi_m\mid\xi_{m+1},\ldots,\xi_{m+n}\}$. From Eq. (\ref{co-1018}), we obtain that
\begin{align*}
&\xi_i\xi_j=-(-1)^{|\xi_i||\xi_j|}\xi_j\xi_i.
\end{align*}
 Moreover,  $\frak{A}(m,n;\underline{t})$ is graded-supercommutative.
 Note that our superalgebra $\frak{A}(m,n;\underline{t})$ is not isomorphic to the  (super)algebra of divided power  mentioned in {\cite[P. 248]{BGLL}},  which is supercommutative.

Suppose $V$ is a module of  Lie superalgebra $\frak{g}$,
and $\frak{A}(V;\underline{t}):=\frak{A}(m,n;\underline{t})$ is {the super-divided power algebra} with  variables $\{v_1,\ldots, v_m\mid v_{m+1},\ldots,
v_{m+n}\}$ which is a fixed basis of $V$. Then $\frak{A}(V;\underline{t})$ has a natural $\frak{g}$-module structure by
defining
\begin{align*}
 g.1=0, \quad
g.u^{(\alpha)}=\sum\limits_{i=1}^{m+n}(-1)^{\alpha_{1}+\cdots+\alpha_{i-1}+|v_i|(\alpha_{m+1}+\cdots+\alpha_{i-1})}(g.v_i)u^{(\alpha-\varepsilon_i)}, \ g\in\frak{g}.\end{align*}

 Now consider the dual module $\frak{g}^{*}$ of the adjoint module of $\frak{g}$. Then
$\frak{A}(\frak{g}^*;\underline{t})$  is a $\frak{g}$-module with respect to the moudue action indicated above. Fix an ordered basis of $\frak{g}$
\begin{align}\label{coho-1505171}
\{x_1,\ldots, x_m\mid x_{m+1},\ldots, x_{m+n}\},\end{align} and write
\begin{align*}\{x_1^*,\ldots, x_m^*\mid x_{m+1}^*,\ldots, x_{m+n}^*\}\end{align*}
for the dual basis.  For simplisity,
 write ${x^*_i}^{(a)}$ for  $u^{(a\varepsilon_i)}$, where $0\leq a\in\mathbb{Z}$, $1\leq i\leq m+n$.

Next let us introduce the differential. Suppose $a_{kl}^i$, $1\leq i,k,l\leq m+n$
are the structure constants of $\frak{g}$ with respect to the basis (\ref{coho-1505171}). Let $\mathrm{d}:\frak{A}(\frak{g}^*;\underline{t})\longrightarrow \frak{A}(\frak{g}^*;\underline{t})$
be the linear operator  defined by
\begin{align*}
&\mathrm{d}(1)=0,\\
&\mathrm{d}x^*_i=\sum_{1\leq k< l\leq m+n}(-1)^{|x_k^*||x_l^*|}a_{kl}^ix^*_kx^*_l-\sum_{m+1\leq k\leq m+n}a_{kk}^i{x^*_k}^{(2)}, \  1\leq i\leq m+n,\\
&\mathrm{d}(u^{(\alpha)})=\sum\limits_{j=1}^{m+n}(-1)^{(1+|x_j^*|)(\alpha_1+\cdots+\alpha_{j-1})}\mathrm{d}(x^*_j){u}^{(\alpha-\varepsilon_j)},\ \alpha \in\frak{A}(m,n;\underline{t}).
\end{align*}
Then  $\mathrm{d}$ is a $\frak{g}$-module homomorphism and
\begin{align*}
& \|\mathrm{d}\|=1; \quad |\mathrm{d}|=\bar{0}; \quad \mathrm{d}^2=0; \\
& \mathrm{d}(u^{(\alpha)}u^{(\beta)})=\mathrm{d}(u^{(\alpha)})u^{(\beta)}+(-1)^{\|u^{(\alpha)}\|}u^{(\alpha)}\mathrm{d}(u^{(\beta)}).
\end{align*}
Recall that the \textit{divided power cohomology} of Lie superalgebra $\frak{g}$,
denoted by $DPH^{\bullet}(\frak{g})$, is defined by the cochain complex $(\frak{A}(\frak{g}^*;\underline{t}),
\mathrm{d})$ (see \cite{BGLL}).

A useful tool to compute the cohomology of Lie superalgebras is  the Hochschild-Serre spectral sequence relative to an ideal $I \triangleleft\frak{g}$
\cite{Musson}:
\begin{align}\label{cohoss}
E_2^{k, s}=H^{k}(\frak{g}/I, H^s(I))\Rightarrow H^{k+s}(\frak{g}).
\end{align}
For Heisenberg Lie superalgebra $\frak{g}=\frak{h}_{2m,n}$ or $\frak{ba}_n$, we are concentrate with the ideal $I=\mathbb{F}z$.
Notice that $\frak{g}/I$ is abelian and the action of $\frak{g}/I$ on $I^*$ is trivial. Thus
\begin{align*}
 H^{k}(\frak{g}/I, H^s(I))=\bigwedge^{k}(\frak{g}/I)^*\bigotimes\bigwedge^s(I)^*.
\end{align*}


\section{Characteristic Zero}\label{nco-h003}
Throughout this section  the ground field $\mathbb{F}$ is assumed to be of characteristic zero.
\subsection{Main Results}\label{co-h001}
Let $\frak{G}_{2m,n}$ be the associative superalgebra with $2m+n$ homogeneous generators
$$\{\xi_1,\ldots,\xi_{2m}\mid \eta_1,\ldots,\eta_n\}$$ 
 and  defining relations with $0\leq i, j\leq 2m,\  1\leq k, s\leq n$,
\begin{align}
  &\xi_i\xi_j+\xi_j\xi_i=0; \quad \xi_i\eta_k+\eta_k\xi_i=0; \quad \eta_k\eta_s-\eta_s\eta_k=0;\label{co-h-140104eq1}\\
  &\sum_{i=1}^m\xi_{i}\xi_{m+i}-{2}^{-1}\sum_{i=1}^n\eta_{i}\eta_i=0.\label{co-h-140104eq2}
\end{align}
 Then $\frak{G}_{2m,n}$ is $\mathbb{Z}$-graded by letting $\|\xi_i\|=\|\eta_k\|=1.$

Let $\frak{G}_{n}$ be the associative superalgebra with $2n$ homogeneous generators
$$\{\xi_1,\ldots,\xi_{n}\mid \eta_1,\ldots,\eta_n\}$$
and  defining relations, for $0\leq i, j\leq n$,
\begin{align}
  &\xi_i\xi_j+\xi_j\xi_i=0; \quad \xi_i\eta_j+\eta_j\xi_i=0; \quad \eta_i\eta_j-\eta_j\eta_i=0;\label{co-h-140104eq3}\\
  &\sum_{k=1}^n\eta_{k}\xi_k=0.\label{co-h-140104eq4}
\end{align}
 Similarly, $\frak{G}_{n}$ is $\mathbb{Z}$-graded by letting $\|\xi_i\|=\|\eta_j\|={1}.$

 From Eqs. (\ref{co-h-140104eq1})--(\ref{co-h-140104eq4}) we know that $\frak{G}_{2m,n}$ and $\frak{G}_{n}$ are just the  quotient algebras of
 super-exterior algebras with variables $\{\xi_1,\ldots,\xi_{2m}\mid \eta_1,\ldots,\eta_n\}$ and $\{\xi_1,\ldots,\xi_{n}\mid \eta_1,\ldots,\eta_n\}$
 modulo
 the  ideals generated by the elements
$\sum_{i=1}^m\xi_{i}\xi_{m+i}-{2}^{-1}\sum_{i=1}^n\eta_{i}\eta_i
$ and $\sum_{k=1}^n\xi_k\eta_{k}$, 
respectively.

Let $\mathcal{V}$  be an infinite-dimensional superspace with  a $\mathbb{Z}_2$-homogeneous basis
 $\{v_i\mid i\in\mathbb{Z}^{+}\}$ with $|v_i|=\bar{i}$.
View $\mathcal{V}$  as  a $\mathbb{Z}$-graded
associative superalgebra with trivial multiplication and a $\mathbb{Z}$-grading structure given by  $\|v_i\|=n+i$.

The superalgbra structures and Betti numbers for the cohomology of Heisenberg Lie superalgebras are formulated as follows:
\begin{theorem}\label{Heisenbergth1}
 As a $\mathbb{Z}$-graded superalgebra,   $H^{\bullet}(\frak{h}_{2m,n})$ is isomorphic to
$\frak{G}_{2m,n}$. In particular,
\[\dim H^k(\frak{h}_{2m,n})=\frak{d}^k(2m,n)-\frak{d}^{k-2}(2m,n).\]
\end{theorem}

\begin{theorem}\label{Heisenbergth2}
As $\mathbb{Z}$-graded superalgebras,
$H^{\bullet}(\frak{ba}_{n})$ is isomorphic to the direct sum of superalgebras $\frak{G}_{n}\bigoplus\mathcal{V}$.
In particular,
\begin{align*}
\dim H^k(\frak{ba}_n)=&\frak{d}^k(n,n)-\frak{d}^{k-2}(n,n)\\
&+ \sum_{i=2}^{\lfloor\frac{k-2}{2}\rfloor}(-1)^{i-1}(\frak{d}^{k-2i}(n,n)-\delta_{k-2i,n})+\sum_{i=1}^{k-1}\delta_{k-i,n}+\delta_{k-2,n}.
 \end{align*}
\end{theorem}
\subsection{Proof}\label{co-h003}

We first consider  the  Heisenberg Lie superalgebra
$\frak{h}_{2m,n}$.

\noindent \textbf{Proof of Theorem \ref{Heisenbergth1}.} 
We use the Hochschild-Serre spectral sequence relative to  the ideal $I=\mathbb{F}z$ [see Eq. (\ref{cohoss})]. For $k\in \mathbb{Z}$, we have
\begin{align*}
&E_2^{k,0}=\bigwedge^k(\frak{h}_{2m,n}/\mathbb{F}z)^*; \quad E_2^{k,1}=\bigwedge^k(\frak{h}_{2m,n}/\mathbb{F}z)^*\bigotimes (\mathbb{F}z)^*;\\
&E_2^{k,s}=0, \quad s\not=0, 1.
\end{align*}
Moreover,
\begin{align*}
E_\infty^{k,s}=E_3^{k,s}, \ s=0\mbox{ or }1; \quad E_\infty^{k,s}=0, \ \mbox{ otherwise}.
\end{align*}
Thus it is sufficient to consider the mapping
\begin{align*}
  \mathrm{d}_2^{k,1}: & \bigwedge^k(\frak{h}_{2m,n}/\mathbb{F}z)^*\bigotimes (\mathbb{F}z)^*\longrightarrow
  \bigwedge^{k+2}(\frak{h}_{2m,n}/\mathbb{F}z)^*\\
  & \mathrm{d}_2^{k,1}(f\otimes z^*)=\mathrm{d}z^*\wedge f, \ f\in \bigwedge^k(\frak{h}_{2m,n}/\mathbb{F}z)^*,
\end{align*}
where
\begin{equation*}
  \mathrm{d}z^*=\sum_{i=1}^mx^*_{i}{\wedge} x^*_{m+i}-2^{-1}\sum_{i=1}^ny^*_{i}{\wedge} y^*_i
\end{equation*}
can be viewed as an element in $\bigwedge^2(\frak{h}_{2m,n}/\mathbb{F}z)^*.$
 From a direct computation, we have $\ker{\mathrm{d}_2^{k,1}}=0$.
Then we obtain that \begin{align*}
                     E_3^{k,0}= \bigwedge^k(\frak{h}_{2m,n}/\mathbb{F}z)^*\Big/
                     \Big(\mathrm{d}z^*\wedge\bigwedge^{k-2}(\frak{h}_{2m,n}/\mathbb{F}z)^*\Big); \quad  E_3^{k,1}=0.
                    \end{align*}
  Thus, $H^k(\frak{h}_{2m,n})=\bigwedge^k(\frak{h}_{2m,n}/\mathbb{F}z)^*\Big/
  \Big(\mathrm{d}z^*\wedge\bigwedge^{k-2}(\frak{h}_{2m,n}/\mathbb{F}z)^*\Big)$.
  The proof is complete.
                    \qed\\

Next we consider  the Heisenberg Lie superalgebra
$\frak{ba}_{n}$.
For an odd element $w^*$, write $({w^*})^l$ for $w^*\wedge\cdots\wedge w^*\in\bigwedge^l\frak{ba}^*_{n}$.
For $\frak{ba}_n$ with $1\leq n\in\mathbb{Z}$, put
\begin{align*}
X^*_n=x^*_1{\wedge}\cdots{\wedge}
x^*_n.
\end{align*}
Following from a standard computation, we obtain that
\begin{align*}
  \mathrm{d}z^*=\sum_{j=1}^nx^*_j{\wedge} y^*_j;
  \quad \mathrm{d}z^*\wedge\mathrm{d}z^*=0,
\end{align*}
and then we can view $ \mathrm{d}z^*$ as an element in $\bigwedge^2(\frak{ba}_n/\mathbb{F}z)^*.$
\begin{lemma}\label{co-h-l3} 
Suppose  $\alpha\in\bigwedge^k(\frak{ba}_n/\mathbb{F}z)^*$, $0\leq k\in\mathbb{Z}$, and  $\alpha\wedge\mathrm{d}z^*=0$. Then
\begin{align*}
\alpha=\beta\wedge\mathrm{d}z^*+\delta_{k,n}\lambda
X^*_n,
\end{align*}
for some $\beta\in {\bigwedge}^{k-2}(\frak{ba}_n/\mathbb{F}z)^*$ and
$\lambda\in\mathbb{F}$.
\end{lemma}

\begin{proof}
The conclusion is trivial for $k=0$ and so we assume that
$k\geq 1$.  To show the conclusion, we use induction on $n$. When $n=1$, we may assume that
\[\alpha=\lambda {y^*_1}^{k-1}{\wedge} x^*_1+\mu {y^*_1}^{k},\ \ \lambda, \mu\in\mathbb{F}.
 \]
 Since $\alpha\wedge \mathrm{d}{z^*}=0$, we have $\mu=0$. Then
 $\alpha=-\lambda {y^*_1}^{k-2}{\wedge}x^*_1\wedge y^*_1+\delta_{k,1}\lambda X^*_1$.
 Now,   we may assume  that $n>1$.
Notice that
\begin{align}\label{coho-150406}
&x^*_n{\wedge}\sum_{j=1}^{n}x^*_j{\wedge} y^*_j=x^*_n{\wedge}\sum_{j=1}^{n-1}x^*_j{\wedge} y^*_j.
\end{align}
Suppose
\begin{align*}
\alpha=\alpha_0+\alpha_1{\wedge} x^*_n+\sum_{i=1}^{k}\gamma_i {\wedge}
{y^*_n}^i+\sum_{i=1}^{k-1}\zeta_i{\wedge} {y^*_n}^i{\wedge} x^*_n,
\end{align*}
 where $\alpha_0\in \bigwedge^{k} (\frak{ba}_{n-1}/\mathbb{F}z)^*$, $\alpha_1\in \bigwedge^{k-1} (\frak{ba}_{n-1}/\mathbb{F}z)^*$,
 $\gamma_i\in \bigwedge^{k-i} (\frak{ba}_{n-1}/\mathbb{F}z)^*$ and $\zeta_i\in \bigwedge^{k-i-1} (\frak{ba}_{n-1}/\mathbb{F}z)^*$.
Note that $|\alpha_0|=|\alpha_1|=|\alpha|$ and
$|\gamma_i|=|\zeta_i|=|\alpha|+\bar{i}$. The condition
$\alpha\wedge \mathrm{d}z^*=0$ implies that
\begin{align}
&\alpha_1{\wedge} x^*_n{\wedge}\sum_{j=1}^{n-1}x^*_j{\wedge} y^*_j=0,\label{co-h-eq2}\\
&\zeta_1{\wedge} y^*_n{\wedge} x^*_n{\wedge}\sum_{j=1}^{n-1}x^*_j{\wedge} y^*_j+
\alpha_0{\wedge}x^*_n{\wedge}
y^*_n=0,\label{co-h-eq3}\\
&\sum_{i=2}^{k-1}\zeta_i{\wedge} {y^*_n}^i{\wedge} x^*_n{\wedge}\sum_{j=1}^{n-1}x^*_j{\wedge} y^*_j+
\sum_{i=1}^{k}\gamma_i{\wedge}{y^*_n}^{i}{\wedge}
x^*_n{\wedge}
y^*_n=0.\label{co-h-140114eq1}
\end{align}
From Eq.  (\ref{co-h-eq2}), we have $\alpha_1{\wedge}\sum_{j=1}^{n-1}x^*_j{\wedge} y^*_j=0$,
which implies that, for $\frak{ba}_{n-1}$,  $\alpha_1\wedge\mathrm{d}z^*=0$.
Applying the inductive hypothesis, we may
assume that
\begin{align}\label{co-h-eq8}
\alpha_1=\beta_1{\wedge}\sum_{j=1}^{n-1}x^*_j{\wedge} y^*_j+\delta_{k-1,n-1}\lambda
X^*_{n-1},
\end{align}
for some  $ \beta_1\in \bigwedge^{k-3} (\frak{ba}_{n-1}/\mathbb{F}z)^*$ and
$\lambda\in\mathbb{F}$.
From  Eq. (\ref{co-h-eq3}), we have
\begin{align}
\zeta_1{\wedge}\sum_{j=1}^{n-1}x^*_j{\wedge} y^*_j=-\alpha_0. \label{co-h-0918-eq1}
\end{align}
From Eq. (\ref{co-h-140114eq1}), 
we have
\begin{align}
&\gamma_{k}=\gamma_{k-1}=0,\label{co-h-eq1}\\
&\sum_{i=2}^{k-1}\zeta_i{\wedge} {y^*_n}^{i-1}{\wedge}\sum_{j=1}^{n-1}x^*_j{\wedge} y^*_j=
-\sum_{i=1}^{k-2}\gamma_i{\wedge}{y^*_n}^{i}.\label{co-h-eq5}
\end{align}
By  Eqs. (\ref{coho-150406}), (\ref{co-h-eq8})--(\ref{co-h-eq5}), we
have
\begin{align*}
\alpha=(\beta_1{\wedge}
x^*_n-\sum_{i=1}^{k-1}\zeta_i\wedge {y_n^*}^{i-1}){\wedge}\sum_{j=1}^nx^*_j{\wedge} y^*_j+\delta_{k,n}\lambda
X^*_n.
\end{align*}
The proof is complete.
\end{proof}


\noindent \textbf{Proof of Theorem \ref{Heisenbergth2}.} From Eq. (\ref{cohoss}) with the ideal $I=\mathbb{F}z$,  we know that for $k\in \mathbb{Z}$, $0\leq
i\leq k$,
\begin{align*}
E_2^{k-i,i}=\bigwedge^{k-i}(\frak{ba}_n/\mathbb{F}z)^*\bigotimes(\mathbb{F}z^*)^i.
\end{align*}
We  consider the mapping
\begin{align*}
  \mathrm{d}_2^{k-i,i}: & \bigwedge^{k-i}(\frak{ba}_n/\mathbb{F}z)^*\bigotimes(\mathbb{F}z^*)^i\longrightarrow
  \bigwedge^{k-i+2}(\frak{ba}_n/\mathbb{F}z)^*\bigotimes(\mathbb{F}z^*)^{i-1}\\
  & \mathrm{d}_2^{k-i,i}(f\otimes (z^*)^i)=if\wedge\mathrm{d}z^*\otimes(z^*)^{i-1}, \ f\in \bigwedge^{k-i}(\frak{ba}_{n}/\mathbb{F}z)^*.
\end{align*}
From Lemma \ref{co-h-l3}, we obtain that
\begin{align*}
&E_3^{k,0}=\bigwedge^{k}(\frak{ba}_n/\mathbb{F}z)^*\Big/
\Big(\bigwedge^{k-2}(\frak{ba}_n/\mathbb{F}z)^*\wedge\mathrm{d}z^*\Big);\\
&E_3^{k-i,i}\cong\delta_{k-i,n}\mathbb{F}X_n^*\otimes(z^*)^i, \ 0<i<k.\\
&E_3^{0,k}=0.
\end{align*}
Moreover,
\begin{align*}
E_4^{k-i,i}=E_3^{k-i,i}, \ 0\leq i\leq k.\\
\end{align*}
Then, we have
\begin{align}\label{coho-1504061435}
H^k(\frak{ba}_n)=\bigwedge^{k}(\frak{ba}_n/\mathbb{F}z)^*\Big/
\Big(\bigwedge^{k-2}(\frak{ba}_n/\mathbb{F}z)^*\wedge\mathrm{d}z^*\Big)\bigoplus\bigoplus_{i=1}^{k-1}\delta_{k-i,n}\mathbb{F}X_n^*\otimes(z^*)^i.
\end{align}
Thus,
$H^{\bullet}(\frak{ba}_{n})\cong \frak{G}_n\bigoplus \mathcal{V}$ as
superalgebras.

To compute the dimension of $H^k(\frak{ba}_n)$, we consider the following mapping: for any $k\in\mathbb{Z}$,
\begin{align*}
\psi_{k}: {\bigwedge}^k (\frak{ba}_{n}/\mathbb{F}z)^*\longrightarrow {\bigwedge}^{k+2}
(\frak{ba}_{n}/\mathbb{F}z)^*, \quad \alpha\mapsto\alpha{\wedge}\mathrm{d}z^*.
\end{align*}
We claim
\[\dim\mathrm{Ker}\psi_{k}=\sum_{i=1}^{\lfloor\frac k
2\rfloor}(-1)^{i-1}\left(\dim {\bigwedge}^{k-2i}
(\frak{ba}_{n}/\mathbb{F}z)^*-\delta_{k-2i,n}\right)+\delta_{k,n}.\]
In fact, we can use induction on $k$ and assume
that $k\geq 1$. When $k=1$, for $\alpha\in\mathrm{Ker}\psi_{k}\subset {\bigwedge}^1 (\frak{ba}_{n}/\mathbb{F}z)^*$, from Lemma \ref{co-h-l3}, we obtain
that \[\alpha=\delta_{1,n}\lambda X^*_1.\] Then
$\dim\mathrm{Ker}\psi_{1}=\delta_{1,n}$. When $k\geq2$, notice
that $X^*_n\in\mathrm{Ker}\psi_{n}$.
From Lemma \ref{co-h-l3},
we obtain that
\begin{align*}
\dim\mathrm{Ker}\psi_{k}&=\dim\mathrm{Im}\psi_{k-2}+\delta_{k,n}\\
&=\dim {\bigwedge}^{k-2} (\frak{ba}_{n}/\mathbb{F}z)^*-\dim\mathrm{Ker}\psi_{k-2}+\delta_{t,n}.
\end{align*}
Then the claim holds.
From     Eq. (\ref{coho-1504061435}) we have
\begin{align*}
\dim H^k(\frak{ba}_n)=&\dim {\bigwedge}^k (\frak{ba}_{n}/\mathbb{F}z)^*-\dim{\bigwedge}^{k-2} (\frak{ba}_{n}/\mathbb{F}z)^*\\
&+\dim\mathrm{Ker}\psi_{k-2}+\sum_{i=1}^{k-1}\delta_{k-i,n}.
\end{align*} The proof is complete.
 \qed\\

\section{Characteristic $p>3$}\label{nco-h004}
Throughout this section the ground field is assumed to be of characteristic $p>3$. Put $\delta_a'=1$ when $a\equiv0 \pmod p$ and $\delta_a'=0$ otherwise.
\subsection{Divided Power Cohomology of $\frak{ba}_{n}$}\label{co-3.1}
In this section, we  will determine the divided power cohomology of  Heisenberg Lie superalgebra $\frak{ba}_n$. Fix  an $n+1$-tuples of positive
integers $ \underline{t}=\left(t_{1},t_{2},\ldots,t_{n+1}\right)$.

Let $\frak{G}_{n}$  be the super-divided power algebra   with variables $$\{\xi_1,\ldots, \xi_n\mid \eta_{1},\ldots, \eta_{n}\}$$ and
  defining relations
\begin{align*}
\sum_{i=1}^n\eta_{i}\xi_i=0.
\end{align*}
Then, $\frak{G}_{n}$ is $\mathbb{Z}$-graded by letting $\|\xi_i\|=\|\eta_j\|={1}.$
 we know that  $\frak{G}_{n}$ is just the  quotient algebra of the super-divided power algebra  with variables  $\{\xi_1,\ldots,\xi_{n}\mid
 \eta_1,\ldots,\eta_n\}$ modulo
 the  ideal generated by the elements
 $\sum_{i=1}^n\xi_i\eta_{i}$.

Let $\mathcal{V}$  be a $2^np^{\sum_{i=1}^nt_i-n}(p^{t_{n+1}}-1)$-dimensional superspace with  a $\mathbb{Z}_2$-homogeneous basis
\begin{align*}
\left\{v_{i_1,\ldots,i_s;j}^{\{a_{i_1},\ldots, a_{i_s}, a_{l_1},\ldots, a_{l_{n-s}}\}}\left|\begin{subarray}{c}
                                                                                              1\leq i_1<\cdots <i_s\leq n,\ 0\leq s\leq n;\\
                                                                                              1\leq j\leq p^{t_{n+1}}-1;\\
                                                                                              0\leq a_{i_k}\leq p^{t_{i_{k}}-1}-1,\ 1\leq k\leq s; \\
                                                                                              1\leq a_{l_k}\leq p^{t_{l_{k}}-1},\ 1\leq k\leq n-s
                                                                                            \end{subarray}\right.
\right\},\end{align*}
where $\left|v_{i_1,\ldots,i_s;j}^{\{a_{i_1},\ldots, a_{i_s}, a_{l_1},\ldots, a_{l_{n-s}}\}}\right|=\overline{a_1+\cdots+ a_n-n+s+j}$.
View $\mathcal{V}$  as a $\mathbb{Z}$-graded
associative superalgebra with the trivial multiplication and the $\mathbb{Z}$-grading structure given by  $\left\|v_{i_1,\ldots,i_s;j}^{\{a_{i_1},\ldots,
a_{i_s}, a_{l_1},\ldots, a_{l_{n-s}}\}}\right\|=(a_1+\cdots+a_n)p-n+2s+j$.

We are in position to describe the graded-supercommutative superalgebra on  the full divided power cohomology $DPH^{\bullet}(\frak{ba}_{n})$.
\begin{theorem}\label{Heisenbergth3}
There exists a graded-supercommutative superalgebra  isomorphism between
$DPH^{\bullet}(\frak{ba}_{n})$ and $\frak{G}_{n}\ltimes\mathcal{V}$ which is the semidirect sum of superalgebras with the following multiplication
\begin{align*}
&{\eta_i}^{(ap)}v_{i_1,\ldots,i_s;j}^{\{a_{i_1},\ldots, a_{i_s}, a_{l_1},\ldots, a_{l_{n-s}}\}}
\\
&=(-1)^{sap}\delta_{i\not\in\{i_1,\ldots,i_s\}}\left(\begin{array}{c}
                                       (a+a_i)p-1\\
                                       ap
                                     \end{array}\right)v_{i_1,\ldots,i_s;j}^{\{a_{i_1},\ldots, a_{i_s}, a_{l_1},\ldots,a+a_i,\ldots,
                                     a_{l_{n-s}}\}}.
\end{align*}
\end{theorem}
To prove this theorem, we make the following preparations.
Let \[\{x^*_1,
\ldots, x^*_{n}\mid y^*_1, \ldots, y^*_{n}; z^*\}\]  be the  dual
basis of the basis $\{x_1,
\ldots, x_{n}\mid y_1, \ldots, y_{n}; z\}$ of $\frak{ba}_{n}$.
Notice that
\begin{align}
&\mathrm{d}z^*=\sum_{i=1}^nx^*_iy^*_i;
  \quad \mathrm{d}z^*\mathrm{d}z^*=0;\nonumber\\
&x^*_n\sum_{i=1}^{n}x^*_iy^*_i=x^*_n\sum_{i=1}^{n-1}x^*_iy^*_i;\label{coho1504221}\\
&{y_n^*}^{(ap-1)}\sum_{i=1}^nx^*_iy^*_i={y_n^*}^{(ap-1)}\sum_{i=1}^{n-1}x^*_iy^*_i,\quad 1\leq a\in\mathbb{Z}.\label{coho1504071}
\end{align}
We can view $\mathrm{d}z^*$ as an element in $\frak{A}_2((\frak{ba}_n/\mathbb{F}z)^*; \underline{t})$.
For $0\leq q\leq n$ and $0\leq k\in\mathbb{Z}$, put
\begin{align*}
V_{q,n}^k=\mathrm{span}_{\mathbb{F}}\Big\{x^*_{i_1}\cdots x^*_{i_q}{y^*_{i_1}}^{(a_{i_1}p)}&\cdots{y^*_{i_q}}^{(a_{i_q}p)}
{y^*_{j_1}}^{(a_{j_1}p-1)}{y^*_{j_{n-q}}}^{(a_{j_{n-q}}p-1)}\Big|\\
&\begin{subarray}{l}
(a_1+\cdots+a_n)p-n+2q=k; \ a_{i_s}\geq 0,\ s\in\overline{1,q};\ a_{j_t}\geq 1,\ t\in\overline{1,n-q};\\
\{j_1,\ldots,j_{n-q}\} \mbox{  \small{is the complement  of} } \{i_1,\ldots,i_q\} \mbox{ \small{in} } \{1,\ldots n\}
\end{subarray}
\Big\}.
\end{align*}
\begin{lemma}\label{coho-150407l1}
Suppose  $\alpha\in\frak{A}_k((\frak{ba}_n/\mathbb{F}z)^*; \underline{t})$, $0\leq k\in\mathbb{Z}$, and  $\alpha\mathrm{d}z^*=0$. Then
\begin{align*}
\alpha=\beta\mathrm{d}z^*+\sum_{i=0}^n\delta'_{k+n-2i}\varepsilon_{i,n}^k
\end{align*}
for some $\beta\in \frak{A}_{k-2}((\frak{ba}_n/\mathbb{F}z)^*;\underline{t})$ and
$\varepsilon_{i,n}^k\in V_{i,n}^k$.
\end{lemma}

\begin{proof}
The conclusion is trivial for $k=0$ and then we assume that
$k\geq 1$.  We use induction on $n$. When $n=1$, we may assume that
$\alpha=\lambda {y^*_1}^{(k-1)} x^*_1+\mu {y^*_1}^{(k)}$ with $\lambda, \mu\in\mathbb{F}$.

 Since $\alpha \mathrm{d}{z^*}=0$, we have
 $(k+1)\mu=0$. We consider the following cases:\\

 \noindent\textbf{Case 1:} $k+1\equiv0 \pmod p$. Note that $k-1\not\equiv0 \pmod p$ for $p>3$. Then
 \begin{align*}
\alpha=\beta\mathrm{d}z^*+\mu {y_1^*}^{(k)}
 \end{align*}
where $\beta=-\frac{\lambda}{k-1}{y_1^*}^{(k-2)}\in\frak{A}_{k-2}((\frak{ba}_1/\mathbb{F}z)^*;\underline{t})$ and ${y_1^*}^{(k)}\in V_{0,1}^k.$

 \noindent\textbf{Case 2:} $k+1\not\equiv0 \pmod p$. Then $\mu=0$. When $k-1\equiv0 \pmod p$, we have
 $$\alpha=\lambda{y_1^*}^{(k-1)}x_1^*\in V_{1,1}^k.$$ When $k-1\not\equiv0 \pmod p$, we have
 $\alpha=\beta\mathrm{d}z^*,$
 where $$\beta=-\frac{\lambda}{k-1}{y_1^*}^{(k-2)}\in\frak{A}_{k-2}((\frak{ba}_1/\mathbb{F}z)^*;\underline{t}).$$
 Thus, the conclusion holds.

 Now,   we may assume  that $n>1$.
Suppose
\begin{align*}
\alpha=\alpha_0+\alpha_1 x^*_n+\sum_{i=1}^{k}\gamma_i
{y^*_n}^{(i)}+\sum_{i=1}^{k-1}\zeta_i {y^*_n}^{(i)} x^*_n,
\end{align*}
 where $\alpha_0\in \frak{A}_{k}((\frak{ba}_{n-1}/\mathbb{F}z)^*;\underline{t})$, $\alpha_1\in \frak{A}_{k-1}((\frak{ba}_{n-1}/\mathbb{F}z)^*;\underline{t})$,
 $\gamma_i\in \frak{A}_{k-i}((\frak{ba}_{n-1}/\mathbb{F}z)^*;\underline{t})$ and $\zeta_i\in \frak{A}_{k-i-1}((\frak{ba}_{n-1}/\mathbb{F}z)^*;\underline{t})$.
Note that $|\alpha_0|=|\alpha_1|=|\alpha|$ and
$|\gamma_i|=|\zeta_i|=|\alpha|+\bar{i}$. The condition
$\alpha{} \mathrm{d}z^*=0$ implies that
\begin{align}
&\alpha_1\sum_{j=1}^{n-1}x^*_j y^*_j=0,\label{coho-1504111}\\
&\gamma_i\sum_{j=1}^{n-1}x^*_j y^*_j=0,\quad 1\leq i\leq k,\label{coho-1504114}\\
&\alpha_0=-\zeta_1\sum_{j=1}^{n-1}x^*_jy^*_j,\label{coho-1504112}\\
&\gamma_i=(-1)^{i+1}(i+1)^{-1}\zeta_{i+1}\sum_{j=1}^{n-1}x^*_jy^*_j,\quad 1\leq i\leq k-1, \ i+1\not\equiv 0\pmod p,\label{coho-1504113}\\
&\zeta_{i+1}\sum_{j=1}^{n-1}x^*_jy^*_j=0,\quad 1\leq i\leq k-2, \ i+1\equiv 0\pmod p,\label{coho-1504115}\\
& \gamma_k=\gamma_{k-1}=0.\label{coho-1504116}
\end{align}
From Eqs.  (\ref{coho-1504111}), (\ref{coho-1504114}), and (\ref{coho-1504115})
we obtain that, for $\frak{A}((\frak{ba}_{n-1}/\mathbb{F}z)^*;\underline{t})$,   $1\leq i\leq k-2,$  $2\leq m\leq k-1, \ i, q\equiv 0\pmod p,$
$$\alpha_1\mathrm{d}z^*=\gamma_i\mathrm{d}z^*=\zeta_{q}\mathrm{d}z^*=0.$$
Applying the inductive hypothesis, we may
assume that
\begin{align}
&\alpha_1=\theta\sum_{j=1}^{n-1}x^*_j y^*_j+\sum_{l=0}^{n-1}\delta'_{k+n-2l-2}\epsilon_{l,n-1}^{k-1},\label{coho-15041110}\\
&\gamma_i=\eta_i\sum_{j=1}^{n-1}x^*_j y^*_j+\sum_{l=0}^{n-1}\delta'_{k+n-2l-i-1}\sigma_{l,n-1}^{k-i},\label{coho-15041111}\\
&\zeta_q=\xi_q\sum_{j=1}^{n-1}x^*_j y^*_j+\sum_{l=0}^{n-1}\delta'_{k+n-2l-q-2}\tau_{l,n-1}^{k-q-1},\label{coho-15041112}
\end{align}
where
\begin{align*}
&\theta\in \frak{A}_{k-3}((\frak{ba}_{n-1}/\mathbb{F}z)^*;\underline{t}),\ \epsilon_{l,n-1}^{k-1}\in V_{l,n-1}^{k-1};\\
&\eta_i\in \frak{A}_{k-i-2}((\frak{ba}_{n-1}/\mathbb{F}z)^*;\underline{t}),\   \sigma_{l,n-1}^{k-i}\in V_{l,n-1}^{k-i};\\
&\xi_q\in\frak{A}_{k-q-3}((\frak{ba}_{n-1}/\mathbb{F}z)^*;\underline{t}),\ \tau_{l,n-1}^{k-q-1}\in V_{l,n-1}^{k-q-1}.
\end{align*}
Note that
\begin{align}
&\overline{\epsilon}_{l+1,n}^k:=\epsilon_{l,n-1}^{k-1}x^*_n\in V_{l+1,n}^k\label{coho-15041113}\\
&\overline{\sigma}_{l,n}^k:=\sigma_{l,n-1}^{k-i+1}{y^*_n}^{(i-1)}\in V_{l,n}^k\label{coho-15041114}\\
&\overline{\tau}_{l+1,n}^k:=\tau_{l,n-1}^{k-q-1}{y^*_n}^{(q)}x^*_n\in V_{l+1,n}^k\label{coho-15041115}
\end{align}
From Eqs. (\ref{coho1504221}), (\ref{coho1504071}), (\ref{coho-1504112})--(\ref{coho-15041115}), we obtain that
$$\alpha=\beta\mathrm{d}z^*+\sum_{i=0}^n\delta'_{k+n-2i}\varepsilon_{i,n}^k$$
 where
\begin{align*}
&\beta=\sum_{i=0
}^{k-2}(\delta'_{i+1}-1)(i+1)^{-1}\zeta_{i+1}{y^*_n}^{(i)}
+\theta x_n^*+\sum_{i=1
}^{k-2}\delta'_{i+1}\eta_{i}{y^*_n}^{(i)}
-\sum_{i=2
}^{k-1}\delta'_i\xi_{i}{y^*_n}^{(i)}x_n^*,\\
&\varepsilon_{i,n}^k=\overline{\epsilon}_{i,n}^k+\overline{\sigma}_{i,n}^k+\overline{\tau}_{i,n}^k.
\end{align*}
\end{proof}
\noindent \textbf{Proof of Theorem \ref{Heisenbergth3}.} From Eq. (\ref{cohoss}) with the ideal $I=\mathbb{F}z$,  we know that for $k\in \mathbb{Z}$, $0\leq
i\leq k$,
\begin{align*}
E_2^{k-i,i}=\frak{A}_{k-i}((\frak{ba}_{n}/\mathbb{F}z)^*;\underline{t})\bigotimes(\mathbb{F}{z^*}^{(i)}).
\end{align*}
We  consider the mapping
\begin{align*}
  \mathrm{d}_2^{k-i,i}: & \frak{A}_{k-i}((\frak{ba}_{n}/\mathbb{F}z)^*;\underline{t})\bigotimes(\mathbb{F}{z^*}^{(i)})\longrightarrow
  \frak{A}_{k-i+2}((\frak{ba}_{n}/\mathbb{F}z)^*;\underline{t})\bigotimes(\mathbb{F}{z^*}^{(i-1)})\\
  & \mathrm{d}_2^{k-i,i}(f\otimes {z^*}^{(i)})=f\mathrm{d}z^*\otimes{z^*}^{(i-1)}, \ f\in \frak{A}_{k-i}((\frak{ba}_{n}/\mathbb{F}z)^*;\underline{t}).
\end{align*}
From Lemma \ref{coho-150407l1}, similar to the proof of Theorem \ref{Heisenbergth2}, we can obtain that
\begin{align*}
DPH^k(\frak{ba}_n)=&\frak{A}_{k}((\frak{ba}_{n}/\mathbb{F}z)^*;\underline{t})\Big/
\frak{A}_{k-2}((\frak{ba}_{n}/\mathbb{F}z)^*;\underline{t})\mathrm{d}z^*\\
&\bigoplus\bigoplus_{\begin{subarray}{c}
                                                                    1\leq i\leq k-1\\
                                                                    0\leq j\leq n
                                                                  \end{subarray}}
\delta'_{k-i+n-2j}V_{j,n}^{k-i}\otimes{z^*}^{(i)}.
\end{align*}
Put
\begin{align*}
&\varepsilon_{i_1,\ldots,i_s;j}^{\{a_{i_1},\ldots, a_{i_s}, a_{l_1},\ldots, a_{l_{n-s}}\}}\\
=&x^*_{i_1}\cdots x^*_{i_s}{y_{i_1}^*}^{(a_{i_1}p)}\cdots {y_{i_s}^*}^{(a_{i_s}p)}{y_{l_1}^*}^{(a_{l_1}p-1)}\cdots
{y_{l_{n-s}}^*}^{(a_{l_{n-s}}p-1)}\otimes {z^*}^{(j)}\in V_{s,n}^{k}\bigotimes\mathbb{F}{z^*}^{(j)},
\end{align*}
where  $\sum_{i=1}^na_ip-n+2s=k$.
Notice that for $\varepsilon_{k_1,\ldots,k_d;q}^{\{b_{k_1},\ldots, b_{k_d},
b_{r_1},\ldots, b_{r_{n-d}}\}}\in V_{d,n}^{t}\bigotimes\mathbb{F}{z^*}^{(q)},$
\begin{align*}
\varepsilon_{i_1,\ldots,i_s;j}^{\{a_{i_1},\ldots, a_{i_s}, a_{l_1},\ldots, a_{l_{n-s}}\}}\varepsilon_{k_1,\ldots,k_d;q}^{\{b_{k_1},\ldots, b_{k_d},
b_{r_1},\ldots, b_{r_{n-d}}\}}\in\frak{A}((\frak{ba}_{n}/\mathbb{F}z);\underline{t})\mathrm{d}z^*\bigotimes\mathbb{F}{z^*}^{(j+q)}.
\end{align*}
For any $i=1,\ldots,n,$ notice that
\begin{align*}
&x_i^*\varepsilon_{i_1,\ldots,i_s;j}^{\{a_{i_1},\ldots, a_{i_s}, a_{l_1},\ldots, a_{l_{n-s}}\}}=0,\quad
y_i^*\varepsilon_{i_1,\ldots,i_s;j}^{\{a_{i_1},\ldots, a_{i_s}, a_{l_1},\ldots, a_{l_{n-s}}\}}=0,
\\
&{y_i^*}^{(ap)}\varepsilon_{i_1,\ldots,i_s;j}^{\{a_{i_1},\ldots, a_{i_s}, a_{l_1},\ldots, a_{l_{n-s}}\}}
\\
&=(-1)^{sap}\delta_{i\not\in\{i_1,\ldots,i_s\}}\left(\begin{array}{c}
                                       (a+a_i)p-1\\
                                       ap
                                     \end{array}\right)\varepsilon_{i_1,\ldots,i_s;j}^{\{a_{i_1},\ldots, a_{i_s}, a_{l_1},\ldots,a+a_i,\ldots,
                                     a_{l_{n-s}}\}}\in V_{s,n}^{k+ap}\bigotimes\mathbb{F}{z^*}^{(j)}.
\end{align*}
Then we have a graded-supercommutative superalgebra isomorphism
$$DPH^{\bullet}(\frak{ba}_{n})\cong \frak{G}_n\ltimes \mathcal{V}.$$
The proof is complete.
 \qed


\subsection{Divided Power Cohomology of $\frak{h}_{2m,n}$ in low dimensions}\label{co-3.0}
For $\frak{h}_{2m,n}$, we only describe the superalgebra structure
of $DPH^\bullet(\frak{h}_{0,1})$ and $DPH^\bullet(\frak{h}_{2,1})$. In this case,  let  $\underline{t}=t$ be a positive integer. Determine the divided power cohomology $DPH^\bullet(\frak{h}_{2m,n})$ with arbitrary $m$ and $n$ should be a more cumbersome and excruciating task. A complete description of $DPH^\bullet(\frak{h}_{2m,n})$    is beyond the scope of this paper.

\begin{example}\label{coho-example1}Recall  $\frak{h}_{0,1}=\mathrm{span}_{\mathbb{F}}\{z\mid y\}$ with $[y,y]=z$.
Note that
\begin{align*}
  &\mathrm{d}y^*=0;\ \mathrm{d}z^*=-{y^*}^{(2)}; \\
  &\mathrm{d}(z^*{y^*}^{(i)})=-\left(
                               \begin{array}{c}
                                 i+2 \\
                                 2 \\
                               \end{array}
                             \right)
  {y^*}^{(i+2)}, \ 1\leq i\leq p^{t}-1.
\end{align*}
Then we obtain that
\begin{align*}
\ker\mathrm{d}=\bigoplus_{i=0}^{p^t-1}\mathbb{F}{y^*}^{(i)}\bigoplus\bigoplus_{i=0}^{p^t-1}(\delta'_{i+1}+\delta'_{i+2})\mathbb{F}z^*{y^*}^{(i)},
\end{align*}
and
\begin{align*}
\mathrm{im}\,\mathrm{d}=\bigoplus_{i=0}^{p^t-1}(1-\delta'_i-\delta'_{i-1})\mathbb{F}{y^*}^{(i)}.
\end{align*}

We construct following surperalgebras in order to describe $DPH^\bullet(\frak{h}_{0,1})$.
Let $\frak{G}_{0,1}$ be a $2p^{t-1}$-dimensional  graded-supercommutative superalgebra with basis
$$\{l_0, m_0,\cdots,l_{p^{t-1}-1}, m_{p^{t-1}-1}\},$$
 where $|l_i|=\bar{i}$, $\|l_i\|=ip$, $|m_i|=\overline{{i+1}}$, $\|m_i\|=ip+1$
and multiplication
\begin{align*}
l_il_j=\left(\begin{array}{c}
                                 (i+j)p \\
                                 ip
                               \end{array}\right)l_{i+j},\quad
l_im_j=\left(
                                        \begin{array}{c}
                                          (i+j)p+1\\
                                          ip \\
                                        \end{array}
                                      \right)m_{i+j}.
\end{align*}
Let $\mathcal{U}$ be a $2p^{t-1}$-dimensional $\mathbb{Z}$-graded superalgebra with $\mathbb{Z}_2$-homogeneous basis
$$\{v_1, u_1, \cdots, v_{p^{t-1}},u_{p^{t-1}}\},$$
 where $|v_i|=\overline{i+1}$, $\|v_i\|=ip$, $|u_i|=\bar{i}$, $\|u_i\|=ip-1$,
and  trivial multiplication.

There exist a graded-supercommutative superalgebra  isomorphism between $DPH^\bullet(\frak{h}_{0,1})$ and $\frak{G}_{0,1}\ltimes\mathcal{U}$
which is the semidirect sum of superalgebras with the multiplication
\begin{align*}
&l_iv_j=\left(
           \begin{array}{c}
             (i+j)p-1 \\
             ip \\
           \end{array}
         \right)v_{i+j};\
          l_iu_j=\left(
           \begin{array}{c}
             (i+j)p-2 \\
             ip \\
           \end{array}
         \right)u_{i+j};\\
&m_iu_j=\left(
           \begin{array}{c}
             (i+j)p-1 \\
             ip+1 \\
           \end{array}
         \right)v_{i+j},
        \end{align*}
where $0\leq i\leq p^{t-1}-1$, $1\leq j\leq p^{t-1}$.
\end{example}

\begin{example}\label{coho-example2}Recall $\frak{h}_{2,1}=\mathrm{span}_{\mathbb{F}}\{z, x_1, x_2\mid y\}$ with $[x_1, x_2]=[y,y]=z$.
Note that
\begin{align*}
  &\mathrm{d}x^*=\mathrm{d}y^*=0;\ \mathrm{d}z^*=x_1^*x_2^*-{y^*}^{(2)}.
\end{align*}
Then we obtain that $\ker\mathrm{d}$ is spanned by the following elements:
\begin{align*}
&x_1^*{y^*}^{(i)},\ x_2^*{y^*}^{(i)},\ x_1^*x_2^*{y^*}^{(i)},\ {y^*}^{(i)},\quad 0\leq i\leq p^t-1;\\
&\delta'_{k+2}\left(\left(
                 \begin{array}{c}
                   k \\
                   2 \\
                 \end{array}
               \right)
{y^*}^{(k)}+x_1^*x_2^*{y^*}^{(k-2)}\right)z^*;\\
&\delta'_{k+1}x_1^*{y^*}^{(k-1)}z^*,\ \delta'_{k+1}x_2^*{y^*}^{(k-1)}z^*,\ \delta'_{k+1}\left(\left(
                 \begin{array}{c}
                   k \\
                   2 \\
                 \end{array}
               \right)
{y^*}^{(k)}+x_1^*x_2^*{y^*}^{(k-2)}\right)z^*;\\
&\delta'_{k}x_1^*{y^*}^{(k-1)}z^*,\ \delta'_{k}x_2^*{y^*}^{(k-1)}z^*,\ \delta'_{k}x_1^*x_2^*{y^*}^{(k-2)}z^*;\ \delta'_{k-1}x_1^*x_2^*{y^*}^{(k-2)}z^*.
\end{align*}
By a direct computation we obtain that
$DPH^\bullet(\frak{h}_{2,1})$ is spanned by
\begin{align*}
&{y^*}^{(ep+s)},\  x_i^*x_j^*{y^*}^{(ep+s)};\\
&x_i^*x_j^*{y^*}^{(kp+s-2)}z^*,\ \left(\left(
                                   \begin{array}{c}
                                     kp+s-2 \\
                                     2 \\
                                   \end{array}
                                 \right){y^*}^{(kp+s-2)}+x_1^*x_2^*{y^*}^{(kp+s-4)}\right)z^*,
\end{align*}
where $x_0:=1$ and $0\leq i<j\leq2,\ 0\leq e\leq p^{t-1}-1, \ 1\leq k\leq p^{t-1}-1,\ s=0$ or 1.
We construct following surperalgebras to describe $DPH^\bullet(\frak{h}_{2,1})$.

Put $\rho^{h,t}_{k,s}=\left(
                        \begin{array}{c}
                          (k+h)p+s+t \\
                          kp+s \\
                        \end{array}
                      \right)
$ and  $\rho^{h,t}=\left(
                        \begin{array}{c}
                          hp+t-2 \\
                          2 \\
                        \end{array}
                      \right)
$, $0\leq k, h\in\mathbb{Z}$, $s, t\in\mathbb{Z}$.
Let $\frak{G}_{2,1}$ be an $8p^{t-1}$-dimensional  graded-supercommutative superalgebra with basis
$$\{l^{k,s}_{(0,0)}, l_{(0,1)}^{k,s}, l_{(1,0)}^{k,s}, l_{(1,1)}^{k,s}\mid 0\leq k\leq p^{t-1}-1,\  s=0 \mbox{ or } 1\},$$
 where $|l^{k,s}_{(i,j)}|=\overline{k+s}$ and  $\|l^{k,s}_{(i,j)}\|=kp+s+i+j$, $(i,j)=(0,0), (0,1), (1,0)$ or $(1,1)$,
with the  multiplication
\begin{align*}
&l_{(i,j)}^{k,s}l_{(0,0)}^{h,t}=(\delta_{s+t, 1}+\delta_{s+t, 0})\rho_{k,s}^{h,t}l_{(i,j)}^{k+h, s+t};\\
&l_{(0,1)}^{k,s}l_{(1,0)}^{h,t}=(-1)^{k+s}(\delta_{s+t, 1}+\delta_{s+t, 0})\rho_{k,s}^{h,t}l_{(1,1)}^{k+h, s+t}.
\end{align*}
Let $\mathcal{W}$ be an $(8p^{t-1}-8)$-dimensional $\mathbb{Z}$-graded superalgebra with  basis
$$\{w_{(0,0)}^{k,s}, w_{(0,1)}^{k,s}, w_{(1,0)}^{k,s}, w_{(1,1)}^{k,s}\mid 1\leq k\leq p^{t-1}-1,\  s=0 \mbox{ or } 1\},$$
where  $|w^{k,s}_{(i,j)}|=\overline{k+s}$ and  $\|w^{k,s}_{(i,j)}\|=kp+s-1+i+j$, $(i,j)=(0,0), (0,1), (1,0)$ or $(1,1)$,
with  trivial multiplication.

There exist a graded-supercommutative  superalgebra  isomorphism between $DPH^\bullet(\frak{h}_{2,1})$ and $\frak{G}_{2,1}\ltimes\mathcal{W}$
which is the semidirect sum of superalgebras with the multiplication
\begin{align*}
&l_{(0,0)}^{k,s}w_{(0,0)}^{h,t}=\rho^{h,t}\rho_{k,s}^{h, t-2}((\delta_{s+t,0}+\delta_{s+t,1})w_{(0,0)}^{k+h, s+t}+\delta_{s+t,2}w_{(1,1)}^{k+h, 0});\\
&l_{(i,j)}^{k,s}w_{(0,0)}^{h,t}=(\delta_{s+t,0}+\delta_{s+t,1})\rho^{h,t}\rho_{k,s}^{h, t-2}w_{(i,j)}^{k+h, s+t};\\
&l_{(0,0)}^{k,s}w_{(i,j)}^{h,t}=(\delta_{s+t,0}+\delta_{s+t,1})
           \rho_{k,s}^{h, t-2}w_{(i,j)}^{k+h, s+t}.
        \end{align*}
where $0\leq k\leq p^{t-1}-1$, $1\leq h\leq p^{t-1}-1$, $s, t=0$ or 1 and $(i,j)=(0,1), (1,0)$ or $(1,1)$.
\end{example}


\begin{thebibliography}{99}



\bibitem[BGLL]{BGLL}
Bouarroudj, S., Grozman, P., Lebedev, A. and Leites, D., 2010, Divided power (co)homology.
Presentations of simple finite dimensional modular Lie superalgebras with
Cartan matrix. \textit{Homology, Homotopy Appl.},   \textbf{12}(1), 237--278.

\bibitem[BGT]{BGT} Bieliavsky, P., Goursac, A. and Tuynman, G., 2012, Deformation quantization for Heisenberg supergroup. \textit{J. Funct. Anal.}, \textbf{263}(3), 549--603.

\bibitem[CJ]{CJ}
Cairns, G. and Jambor, S., 2008, The cohomology of the Heisenberg
Lie algebras over fields of finite characteristic. \textit{Proc.
Amer. Math. Soc.}, { \bf 136}, 3803--3807.

\bibitem[CKN]{CKN} Chang-Young, E., Kim, H. and Nakajima, H., 2007, Noncommutative superspace and super Heisenberg group. \textit{J. High Energy of Phys.}, \textbf{4}(4), 403--410.


\bibitem[F]{F}
Fuks (Fuchs), D., 1986,  Cohomology of Infinite-Dimensional Lie
Algebras. {Consultants Bureau}, New York.






\bibitem[G1]{G1}
Gruson, C., 1997, Finitude de l'homologie de certains modules de dimension
finie sur une super alg\'ebre de Lie.  (French) [Finiteness of the
homology of certain finite-dimensional modules over a Lie
superalgebra] \textit{Ann.  Inst.  Fourier (Grenoble)}, \textbf{47}(2),
531--553.

\bibitem[G2]{G2}
Gruson, C., 2000, Sur la cohomologie des super alg\`ebres de Lie
\'etranges.  (French) [Cohomology of strange Lie superalgebras]
\textit{Transform.  Groups}, \textbf{5}(1), 73--84.

\bibitem[GL]{GL}
Grozman, P. and  Leites, D., 2004,  Lie superalgebra structures in
$H^{\bullet}(\frak{g};\frak{g})$. \textit{Czech. J. Phys}.,
\textbf{54}, 1313--1319.


\bibitem[GM]{GM}   Goursac, A. and Michel, J., Superunitary representations of Heisenberg supergroups. arXiv: 1601.07387v1 [math.RT] 27 Jan 2016.


\bibitem[Leb1]{Leb1}
Lebedev, A., 2007, On the Bott-Borel-Weil and Tolpygo theorems. (Russian) \textit{Mat. Zametki}, \textbf{81}(3),  474--477; \textit{translation in Math.
Notes}, \textbf{81}(3-4),  417--421.
\bibitem[Leb2]{Leb2}
Lebedev, A., 2010, Analogs of the orthogonal, Hamiltonian, Poisson, and contact Lie superalgebras in characteristic 2. \textit{J. Nonlinear Math. Phys.},
\textbf{17},  217--251.




\bibitem[LLS]{LLS}
Lebedev, A.,  Leites, D. and Shereshevskii, I., 2005, Lie
superalgebra structures in cohomology spaces of Lie algebras with coefficients in the adjoint representation. in: \textit{Lie Groups and
Invariant Theory}, \'{E}. Vinberg (ed.), \textit{Amer. Math. Soc.
Transl. Ser.}, \textbf{2}(213) (\textit{Amer. Math. Soc.},
Providence, RI), 157--172.

\bibitem[Lei]{Leites}
Leites, D., 1975, Cohomology of Lie superalgebras. \textit{Funkt.
Anal. Pril.}, \textbf{9}(4), 75--76.

\bibitem[M]{Musson}
Musson, I., 2012, Lie Superalgebras and Enveloping Algebras. American Mathematical Society, Providence, RI.

\bibitem[RSS]{RSS}
Rodr\'{\i}guez-Vallarte, M., Salgado, G. and
S\'{a}nchez-Valenzuela, O., 2011, Heisenberg Lie superalgebras and
their invariant superorthogonal and supersymplectic forms.
\textit{J. Algebra}, \textbf{332}, 71--86.




\bibitem[Sa]{LJS}
Santharoubane, L., 1983,  Cohomology of Heisenberg Lie algebras.
\textit{Proc. Amer. Math. Soc.}, { \bf 87}(1), 23--28.

\bibitem[Sk]{Emil}
Sk\"oldberg, E., 2005, The homology of Heisenberg Lie algebras over
fields of characteristic two. \textit{Math. Proc. R. Ir. Acad.},
\textbf{105A}(2), 47--49.

\bibitem[T]{T} Tuynman, G., 2010, Super symplectic geometry and prequantization. \textit{J. Geom. Phys.}, \textbf{60}, 1919--1939.


\bibitem[WLS]{WLS}  Linch  III, W. and  Randall, S., 2015, Superspace de Rham complex and relative
cohomology. \textit{J. High Energy Phys.} \textbf{09}, 190.






\end{thebibliography}
\end{document}